\documentclass{notices}
\RequirePackage[final]{microtype}
\usepackage{amsfonts,amssymb,amsmath,amscd,amsthm}
\usepackage{xcolor,mathtools,xspace}
\usepackage{paralist}
\usepackage{enumitem}
\usepackage{hyperref}

\mathtoolsset{showonlyrefs=true,showmanualtags=true}

\definecolor{link1}{rgb}{0,0,.7}
\definecolor{link2}{rgb}{0,0.25,0.5}
\hypersetup{final,colorlinks,allcolors=link1,citecolor=link2}

\newcommand{\lap}{\Delta}

\newcommand{\grad}{\nabla}

\newcommand{\dv}{\grad \cdot}

\renewcommand{\epsilon}{\varepsilon}
\renewcommand{\leq}{\leqslant}
\renewcommand{\geq}{\geqslant}
\renewcommand{\le}{\leqslant}

\newcommand{\R}{\mathbb{R}}

\newcommand{\Z}{\mathbb{Z}}

\newcommand{\T}{\mathbb{T}}

\makeatletter
\newcommand{\GI@given}[1]{\nonscript\:\mathopen{}#1\vert\nonscript\:\mathopen{}}
\newcommand{\given}[1][]{\GI@given{#1}}
\DeclarePairedDelimiterX\paren[1]{(}{)}{%
  \renewcommand{\given}{\GI@given{\delimsize}}#1%
}
\DeclarePairedDelimiterX\brak[1]{[}{]}{%
  \renewcommand{\given}{\GI@given{\delimsize}}#1%
}

\newcommand{\GI@st}[1]{\nonscript\:#1\vert\nonscript\:\mathopen{}\allowbreak}
\newcommand{\st}[1][]{\GI@st{#1}}
\DeclarePairedDelimiterX\set[1]\{\}{%
  \renewcommand{\st}{\GI@st{\delimsize}}#1%
}
\makeatother

\DeclarePairedDelimiter{\abs}{\lvert}{\rvert}
\DeclarePairedDelimiter{\norm}{\lVert}{\rVert}

\newcommand{\defeq}{\stackrel{\scriptscriptstyle\textup{def}}{=}}

\newcommand{\etal}{et\penalty50\ al.\xspace}

\newcommand{\tdis}{t_{\mathrm{dis}}}

\newcommand\de{{\partial}}
\newcommand{\NN}{\mathbb{N}}
\newcommand{\ZZ}{\mathbb{Z}}
\newcommand\TT {{\mathbb T}}
\newcommand\RR {{\mathbb R}}
\newcommand\e{{\rm e}}
\newcommand\dd{{\rm d}}

\newtheorem{proposition}{Proposition}[section]
\newtheorem{theorem}[proposition]{Theorem}

\theoremstyle{definition}

\newtheorem{question}[proposition]{Question}

\numberwithin{equation}{section}

\iftrue
  \onecolumn
  \RequirePackage{geometry}
\fi
\begin{document}
\title{
Mixing in incompressible flows: \\ transport, dissipation, and their interplay.%
\thanks{This work has been partially supported by the US National Science Foundation under grant DMS-2108080  and the CMU Center for Non
Linear Analysis to GI, the US National Science Foundation under grants DMS-1615457, DMS-1909103, DMS-2206453 to ALM, the European Research Council under grant 676675 FLIRT and the Swiss National Science Foundation under grant 212573 FLUTURA to GC, the Royal Society URF\textbackslash R1\textbackslash 191492 and the EPSRC Horizon Europe Guarantee EP/X020886/1 to MCZ.
}
}

\author{
 Michele Coti Zelati
  \affil{Michele Coti Zelati is Reader in the Department of Mathematics at Imperial College, London, UK. His email address is m.coti-zelati@imperial.ac.uk.
    }
  \and
 Gianluca Crippa
  \affil{
    Gianluca Crippa is Associate Professor in Analysis in the Department of Mathematics and Computer Science at University of Basel, Switzerland. His email address is gianluca.crippa@unibas.ch.
   }
   \and
 Gautam Iyer
  \affil{
    Gautam Iyer is Professor of Mathematics at Carnegie Mellon University, Pittsburgh, USA. His email address is gautam@math.cmu.edu.
   }
   \and
Anna L. Mazzucato
  \affil{
   Anna Mazzucato is Professor of Mathematics at Penn State University, University Park, USA. Her email address is alm24@psu.edu.
   }
}

\maketitle

\section{Introduction}\label{s:intro}

Mixing in fluid flows is an ubiquitous phenomenon, and arises in many situations ranging from everyday occurrences, such as mixing of cream in coffee, to fundamental physical processes such as circulation in the oceans and the atmosphere.
From a theoretical point of view, mixing has been studied since the late nineteenth century in different 
contexts, including dynamical systems, homogenization, control, hydrodynamic stability and turbulence theory.
Although certain aspects of these theories still elude us, significant progress has allowed to provide a rigorous mathematical description of some fundamental mixing mechanisms.
In this survey, we address mixing from the point of view of partial differential equations, motivated by applications that arise in fluid dynamics.
A prototypical example is the movement of small tracer particles (e.g.~pollen grains) in a liquid.
One can visually see the particles ``mix'', and our interest is to quantify this phenomenon mathematically and formulate rigorous results in this context.

When the diffusive effects are negligible, the evolution of the density of tracer particles is governed by the transport equation
\begin{equation}\label{e:transport}
  \de_t\rho +u\cdot\nabla \rho=0 \,.
\end{equation}
Here $\rho = \rho(t, x)$ is a scalar representing the density of tracer particles, and~$u = u(t, x)$ denotes the velocity of the ambient fluid.
We will always assume that the ambient fluid is \emph{incompressible}, which mathematically translates to the requirement that~$u$ is divergence free.
Moreover, we will only study situations where~$\rho$ is a \emph{passive scalar} (or passively advected scalar) -- that is, the effect of tracer particles on the advecting fluid is negligible and the evolution of~$\rho$ does not influence the fluid velocity field~$u$.
One example where passive advection arises in nature is when light, chemically non-reactant particles are carried by a large fluid body (e.g.\ plankton blooms in the ocean).
Examples of \emph{active scalars} (i.e.,\ scalars which are not passive) are quantities such as salinity and temperature in geophysical contexts.

Our interest is to study mixing away from boundaries, and hence we will study~\eqref{e:transport} with periodic boundary conditions.
For simplicity, and clarity of exposition, we fix the dimension $d=2$.
We mention, however, that most
of the results we state can be extended to higher dimensions without too much difficulty.

We complement \eqref{e:transport} with an initial condition $\rho_0$ at time $t=0$.
Since~$u$ is divergence free, integrating~\eqref{e:transport} in space shows that the total mass
\begin{equation}
  \bar \rho = \int_{\T^2} \rho(t, x) \, \dd x
\end{equation}
is a conserved quantity (i.e.\ remains constant in time).
Thus replacing~$\rho_0$ with~$\rho_0 - \bar \rho$ if necessary, there is no loss of generality in assuming $\rho_0$ (and hence~$\rho(t, \cdot)$) are mean zero.
We will subsequently always assume that~$\rho$ is mean zero.
Physically, $\rho$ now represents the deviation from the mean of the density of tracer particles.

\begin{figure}[hbt]
  \hfil
  \includegraphics[width=.45\linewidth]{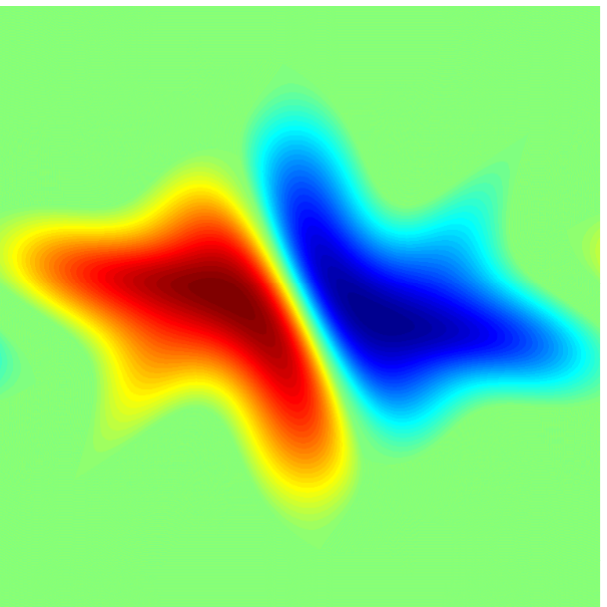}%
  \hfil
  \includegraphics[width=.45\linewidth]{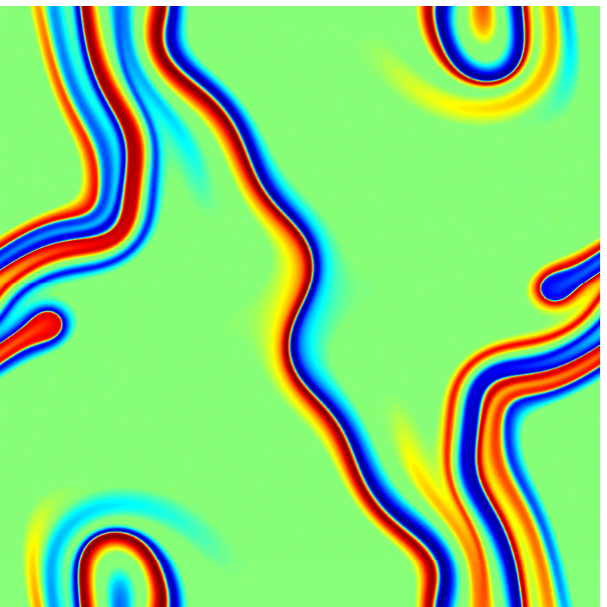}%
  \hfil
  \caption{Example of mixing.
    The red and blue level sets in both figures have exactly the same area.
    For the left figure, averages on scales comparable to $1/8^\text{th}$ of the period are of order~$1$.
    On the right, however, the sets are stretched and interspersed in such a manner that averages at the same scale are much smaller.
  }
  \label{f:mixing}
\end{figure}

Mixing, informally speaking, is the process by which uneven initial configurations transform into a spatially uniform one.
In our setting (since the velocity~$u$ is divergence free), the area of regions of relatively higher (or lower) concentration is preserved.
That is, for any~$c \in \R$, the area of the sub-level sets~$\set{ x \st \rho(t, x) < c }$ and the super-level sets~$\set{x \st \rho(t, x) > c}$ are both constant in time.
Thus if initially the set where~$\rho_0$ is positive occupies half the torus, then for all time the set where~$\rho(t)$ is positive must also occupy half the torus.
The process of mixing will transform~$\rho$ in such a way that the set~$\set{ \rho(t) > 0 }$ will be stretched into many long thin filaments (that still occupy a total area of half), and are interspersed with filaments of the set~$\set{\rho(t) < 0}$ in such a manner that averages at any fixed scale become small (see Figure~\ref{f:mixing}, below).
Mathematically, this is essentially \emph{weak convergence} of~$\rho(t)$ to $0$.
(Recall by assumption~$\bar \rho = 0$.
  If this was not the case, we would instead have weak convergence of~$\rho(t)$ to $\bar \rho$.)

More precisely, we say~$\rho(t)$ becomes mixed as~$t \to \infty$ if~$\rho(t)$ converges weakly to~$\bar \rho$ in~$L^2$.
That is, for every~$L^2$ test function~$f$ we have
\begin{equation}
  \lim_{t \to \infty} \int_{\T^2} \rho(t, x) f(x) \, \dd x = \bar \rho \int_{\T^2} f(x) \, \dd x = 0 \,,
\end{equation}
where the last equality follows because~$\bar \rho = 0$ by assumption.
The standard notation for this convergence is to write
\begin{equation}
  \rho(t) \xrightharpoonup[L^2]{t \to \infty} 0\,.
\end{equation}
As the name (and notation) suggest, this is weaker than \emph{strong~$L^2$ convergence}
\begin{equation}
  \rho(t) \xrightarrow[L^2]{t \to \infty} 0\,,
\end{equation}
which means $\norm{\rho(t)}_{L^2} \to 0$ as~$t \to \infty$.
In our situation given that $u$ is divergence free, we note that for all~$t \in \R$ we have~$\norm{\rho(t)}_{L^2} = \norm{\rho_0}_{L^2}$.
Thus while many of our examples exhibit mixing (i.e.\ weak convergence of~$\rho(t)$ to $0$), they will not have strong convergence of~$\rho(t)$ to $0$, unless~$\rho_0$ is identically~$0$.

Even though weak convergence is a natural way to study
mixing, the disadvantage is that it does not apriori give a quantifiable rate.
To explain further, if~$\rho(t)$ converges to $0$ strongly in~$L^2$, then at time~$t$ the quantity~$\norm{\rho(t)}_{L^2}$ is a measure of how close~$\rho(t)$ is to its (strong) limit.
If~$\rho(t)$ converges to $0$ \emph{weakly}, then~$\norm{\rho(t)}_{L^2}$ may not contain any useful information about the convergence.
(Indeed, for solutions to~\eqref{e:transport},~$\norm{\rho(t)}_{L^2}$ is independent of~$t$.)
It turns out, however, that
in our situation,
weak~$L^2$ convergence to~$0$ is \emph{equivalent} to strong convergence to~$0$ in any \emph{negative Sobolev space}.
Now the norm in these negative Sobolev spaces (which we will define shortly) can be used as a measure of how ``well mixed'' the distribution is (see for instance~\cite{Thiffeault12}).
 
It is easiest to define the negative Sobolev norms using the Fourier series.
Given an (integrable) function~$\theta$ on the torus, we define its Fourier coefficients by
\begin{equation}
  \hat \theta_k = \int_{\T^2} \theta(x) \e^{-2 \pi i \langle x, k \rangle } \, \dd x \,,
  \quad
  \text{where}
  \quad
  k \in \Z^2 \,.
\end{equation}
For mean-zero functions the $0$-th Fourier coefficient,~$\hat \theta_0$, vanishes.
Now, for any~$s \in \R$  we define the homogeneous Sobolev norm of index~$s$ by
\begin{equation} \label{e:homSobolev}
  \|\theta\|_{\Dot{H}^s}^2 \defeq \sum_{k \in \mathbb{Z}^2, k\ne 0} |k|^{2s} |\Hat{\theta}_k|^2 \,.
\end{equation}
Note that for~$s > 0$ the norm puts more weight on higher frequencies.
Thus functions that have a smaller fraction of their Fourier mass in the high frequencies will be ``less oscillatory'' and have a smaller~$\dot H^s$ norm.
For~$s < 0$, however, the norm puts less weight on higher frequencies.
Thus functions that have a larger fraction of their Fourier mass in the high frequencies will be ``very oscillatory'', and have a smaller~$\dot H^s$ norm.
This is consistent with what we expect from ``mixed'' distributions.
Moreover, the result mentioned previously guarantees that mixing of~$\rho$ is equivalent to
\begin{equation}
  \norm{\rho(t)}_{\dot H^s} \xrightarrow{t \to \infty} 0\,,
  \quad\text{for every } s < 0\,,
\end{equation}
Thus, for any~$s < 0$,  the quantity~$\norm{\rho(t)}_{\dot H^s}$ can be used as a measure of how ``mixed'' the distribution is at time~$t$.

For this reason negative Sobolev norms are often referred to as ``mix norms''.
Choosing $s=-1$ is particularly convenient, as the ratio of the $\dot H^{-1}$ norm to the $L^2$ norm scales like a length, and typically represents a length scale of the largest unmixed region.
Since the $L^2$ norm is preserved by equation~\eqref{e:transport}, we can identify the $\dot H^{-1}$ norm with the {\em mixing scale} of the scalar field $\rho$.
We mention that there is also a related notion of mixing scale, which is more geometric in nature (see for instance~\cite{Bressan03}),
but when studying evolution equations the mix-norms described above are easier to work with.

Practically, in order to mix a given initial configuration to a certain degree, one has to expend energy by ``stirring the fluid''.
A natural, physically meaningful question, is to bound the mixing efficiency~\cite{LinThiffeaultEA11}.
That is, given a certain ``cost'' associated with stirring the ambient fluid, what is the most efficient way to mix a given initial configuration?
Two cost functions that are particularly interesting are the energy~$\norm{u}_{L^2}^2$ (which is proportional to the actual \emph{kinetic energy} of the fluid, assuming the fluid is homogeneous), or the enstrophy~$\norm{u}_{\Dot{H}^1}^2$ (which is proportional to the 
fluid's viscous power dissipation).
Hence our question about mixing efficiency can now be formulated as follows.

\begin{question}\label{q:constrained}
  What are optimal bounds on $\norm{\rho(t)}_{\Dot{H}^{-1}}$ in terms of the fluid's energy or enstrophy?
\end{question}


Before answering this question, we note that in order to guarantee solutions to~\eqref{e:transport} exist in the classical sense (i.e.,\ $\rho \in C^{1,1}$ and satisfies~\eqref{e:transport} at every point), we need the fluid velocity field~$u$ to be Lipschitz.
However, the natural first constraints on~$u$ described above do not require~$u$ to be Lipschitz.
Moreover, one intuitively expects efficient mixing flows to be ``turbulent'', and standard turbulence models have velocity fields that are only H\"older continuous at every point.
Thus, in many natural situations arising in the study of fluids, one has to study~\eqref{e:transport} when the advecting velocity field is not Lipschitz.
Seminal work of DiPerna and Lions in '89 addresses this, and shows that certain ``renormalized'' solutions to~\eqref{e:transport} are unique, provided~$u, \grad u \in L^1$.
(This was later extended to a larger class of functions by Ambrosio and we refer the reader to~\cite{Ambrosio04} and references therein for details.)
\smallskip

Returning to Question~\ref{q:constrained}, one can use direct energy estimates to show that if the fluid is \emph{energy constrained} (i.e.\ if~$\norm{u(t)}_{L^2}^2 \leq E$, for some constant~$E$), then~$\norm{\rho(t)}_{H^{-1}}$ can decrease at most linearly as a function of time.
An elegant slice and dice construction of Bressan~\cite{Bressan03} provides an example where this bound is indeed attained (see~\cite{LunasinLinEA12}).
In particular, this provides an example where for some~$T < \infty$ and an incompressible, finite-energy velocity field~$u$,  we have~$\rho(t) \to 0$ weakly in~$L^2$ as~$t \to T$.
That is, the fluid mixes the initial configuration ``perfectly'' in finite time.

On the other hand, if one imposes an~\emph{enstrophy constraint} (i.e.\ a restriction on the growth of~$\norm{u(t)}_{H^1}^2$), or more generally a restriction on the growth of~$\norm{u}_{L^p} + \norm{\grad u(t)}_{L^p}$, then the DiPerna--Lions theory guarantees finite-time perfect mixing can not occur.
Indeed, equation~\eqref{e:transport} is time reversible, and so 
finite-time perfect mixing would provide one non-trivial solution to~\eqref{e:transport} with initial data~$0$.
Since~$\rho \equiv 0$ is clearly another solution, we have non-uniqueness for weak solutions to~\eqref{e:transport}, which is not allowed by the  DiPerna-Lions theory when~$u, \grad u \in L^1$. So one can not have finite-time perfect mixing in this case.

More quantitatively, one can use the regularity of DiPerna--Lions flows~\cite{CrippaDeLellis08} to obtain explicit exponential lower bounds on the mix norm.
Namely, one can prove~\cite{IyerKiselevEA14}
\begin{equation}\label{e:mixlower}
  \norm{\rho(t)}_{H^{-1}}
    \geq C_0 \exp\paren[\Big]{ - C_1 \int_0^t \norm{\grad u(s)}_{L^p} \, \dd s } \,,
\end{equation}
for every~$p \in (1, \infty]$, and some constants $C_0, C_1$ that depend on~$\rho_0$ and~$p$.
(We remark that, for~$p = \infty$, an elementary proof of the lower bound~\eqref{e:mixlower} follows from  Gronwall's inequality and the method of characteristics.
For~$p \in (1, \infty)$, however, the proof is more involved and requires some tools from geometric measure theory.)


Interestingly, whether or not~\eqref{e:mixlower} holds for~$p = 1$ is an open question.
Indeed, the proof of the needed regularity estimates for the flow in~\cite{CrippaDeLellis08} relies on boundedness of a maximal function, which fails for~$p =1$.
The bound \eqref{e:mixlower} for $p=1$ is related to a conjecture of Bressan~\cite{Bressan03} on the cost of rearranging a set, which is still an open question.

For optimality, there are now several constructions of velocity fields that show~\eqref{e:mixlower} is sharp.
These constructions produce enstrophy constrained velocity fields for which~$\norm{\rho(t)}_{\Dot{H}^{-1}}$ decays exponentially in time.
A construction in~\cite{AlbertiCrippaEA19} does this by starting with initial data which is supported in a strip and finds a Lipschitz velocity field that pushes it along a space filling curve.
Constructions in~\cites{BCZG22,BedrossianBlumenthalEA21,HillSturmanEA22} produce regular velocity fields for which
\begin{equation}\label{e:expmix}
  \norm{\rho(t)}_{\dot H^{-1}} \leq D\, \e^{-\gamma t} \norm{\rho_0}_{\dot H^1}\,,
\end{equation}
for every initial data~$\rho_0 \in \dot H^1$.
Such flows are called \emph{exponentially mixing}, and we revisit this in more detail in Section~\ref{s:om}.
\medskip

In addition to optimal mixing, there are three other themes discussed in this article.
We briefly introduce these themes here, and elaborate on them in subsequent sections.
\begin{asparaenum}
  \item \emph{Loss of regularity.}
    When~$u$ is regular, classical theory guarantees regularity of the initial data~$\rho_0$ is propagated by the equation~\eqref{e:transport}.
    However, when~$u$ is irregular (e.g.\ when~$\grad u \in L^p$ with~$p < \infty$), it may happen that all regularity of the initial data~$\rho_0$ is immediately lost.
    Not surprisingly, this loss is intrinsically related to mixing.
    Indeed, the process of mixing generates high frequencies, making the initial data more irregular.
    When~$u$ is not Lipschitz one can arrange rapid enough growth of high frequencies to ensure that \emph{all} Sobolev regularity of the initial data is immediately lost.
    We describe this construction in detail in Section~\ref{s:reg}.

  \item \emph{Enhanced dissipation.}
    In several physically relevant situations, both diffusion and transport are simultaneously present. The nature of diffusion is to~rapidly damp high frequencies. Since mixing generates high frequencies, the combined effect of mixing and diffusion will lead to energy decay of solutions that is an order of magnitude faster than when diffusion acts alone.
    This phenomenon is known as~\emph{enhanced dissipation} and is described in Section~\ref{s:ed}.

  \item \emph{Anomalous dissipation.}
    Even under the enhanced dissipation mentioned above, for which  the energy decay is much faster due to the combined effect of diffusion and transport, the energy decay rate vanishes with the diffusivity.
    In some sense, this outcome is expected, as solutions to~\eqref{e:transport} (formally) conserve energy.
    For certain (irregular) flows, however, it is possible for the energy decay rate in the presence of small diffusivity to stay uniformly positive,
    a phenomenon known as~\emph{anomalous dissipation}. It implies, in particular,
    that the vanishing-diffusivity limit can produce
    \emph{dissipative solutions} of~\eqref{e:transport},
    that is, solutions for which the energy decreases with time.
    We discuss anomalous dissipation in Section~\ref{s:ad}.
\end{asparaenum}

\section{Optimally mixing flows} \label{s:om}
In this section, our primary focus centers around understanding the concept of \emph{shearing} as one of the central mechanisms of mixing, and how this mechanism gives rise to flows that mix optimally.

\subsection{Shear flows}
Shear flows are the simplest example of incompressible flows on $\TT^2$. Their streamlines (lines tangent to the direction of the velocity vector) are parallel to each other and the velocity takes the form $u=(v(x_2),0)$. The corresponding transport equation is
\begin{equation}\label{eq:sheartransport}
\de_t \rho+v(x_2)\de_{x_1} \rho=0,\qquad \rho(0,x)=\rho_0(x),
\end{equation}
the solution $\rho(t,x_1,x_2)=\rho_0(x_1-v(x_2)t,x_2)$ of which can be computed explicitly via the method of characteristics.

If the initial datum only
depends on $x_2$, then the solution remains constant for all times. Otherwise, a hint of creation of small scales is given by 
the growth of $\|\de_{x_2}\rho\|_{L^2}$ linearly in time. To deduce a quantitative mixing estimate, one can take a partial Fourier 
transform in $x_1$ of \eqref{eq:sheartransport}: denoting $\hat{\rho}(t,k,x_2)$, with $k\in \ZZ$, the Fourier coefficients of $\rho$,  
\eqref{eq:sheartransport} becomes
\begin{equation}\label{eq:sheartransportFourier}
\de_t \hat{\rho}+ikv(x_2) \hat{\rho}=0,\qquad \hat{\rho}(0,k,x_2)=\hat{\rho}_0(k,x_2).
\end{equation}
Since $ \hat{\rho}(t,k,x_2)=\e^{-ikv(x_2)t}\hat{\rho}_0(k,x_2)$, mixing follows by estimating oscillatory integrals of the form
\begin{equation}\label{eq:corrdecay}
\int_{\TT}\e^{-ikv(x_2)t}\hat{\rho}_0(k,x_2) \hat{\phi}(k,x_2)\dd x_2, \qquad \phi \in \dot{H}^1(\TT^2).
\end{equation}
A duality argument and an application of the stationary-phase lemma entails a (sharp) mixing estimate.

\begin{theorem}\label{thm:mixingshear}
Assume that $v \in C^m(\TT)$, for some integer $m\geq 2$, and its derivatives up to order $m$ do not vanish 
simultaneously:
$|v'(x_2)| + |v''(x_2)| + \dots + |v^{(m)}(x_2)| > 0 $, for all  $x_2 \in \TT$. 
Then there exists a positive constant $C=C(v)$ such that
\begin{equation}
  \|\rho(t)\|_{\dot{H}^{-1}} \le \frac{C}{t^{1/m}}\|\rho_0\|_{\dot{H}^1}, \qquad \forall t\geq 1,
\end{equation}
for all initial data $\rho_0 \in \dot{H}^1(\TT^2)$ with vanishing $x_1$-average. 
\end{theorem}

The mixing rate is solely determined by how degenerate the critical points of $v$ are, and Theorem~\ref{thm:mixingshear} 
tells us that the flatter the critical points, the slower the (universal) mixing rate. 
In the case of simple critical points, such as the Kolmogorov flow
$v(x_2)=\sin x_2$, we have $m=2$. 

The requirement that $\rho_0$ has vanishing $x_1$-average (namely, $k\neq 0$ in  \eqref{eq:sheartransportFourier})
is essential: it excludes
functions that are constant on streamlines, i.e., eigenfunctions (with eigenvalue 0) of the transport operator, 
which do not enjoy any mixing. 

\subsection{General two-dimensional flows}
Although regular shear flows can achieve algebraic mixing rates, we could be inclined to think that their simple structure constitutes an obstruction to faster mixing. It turns out that if $H$ is an autonomous, non-constant Hamiltonian function on $\TT^2$, of class $C^2$,
generating an
incompressible velocity field $u=\nabla^\perp H =(-\de_{x_2}H,\de_{x_1}H)$, then the mixing rate of $u$ is at best $1/t$, and can be even slower depending on the structure of $H$, see \cite{BCZM22}. This result can be interpreted as follows:
despite the fact that $H$ could have hyperbolic points, at which the flow map displays exponential stretching and compression, shearing is the main mixing
mechanism in 2d. This can be deduced from the existence of an invariant domain for the Lagrangian flow on which $u$  is bounded away from zero, 
which in turn implies that there exists a well-defined, regular and invertible change of coordinates  $(x_1,x_2)\mapsto (h,\theta)\in \TT\times (h_0,h_1)$,
where the interval $(h_0,h_1)\subset {\rm range}(H)$ is determined by the invariant set.

  \begin{proposition}\label{prop:keystep}
 Let $H\in C^2(\TT^2)$.
  There exists an invariant  open set $\mathcal{I} 	\subset \TT^2$ such that for any $\rho_0 \in C^1(\TT^2)$ with 
	${\rm supp} (\rho_0) \subset \mathcal{I}$, the corresponding solution $\rho$ of \eqref{e:transport} satisfies  
  	\begin{equation}\label{eq:reg}
  		\|  \rho( t) \|_{\dot{H}^1}
  		\le C (1+t)\|\nabla \rho_0\|_{L^\infty},
  	\end{equation}
for some $C=C(\mathcal{I}, H)$ and all $t\geq 0$.
  \end{proposition}

The set of coordinates $(h,\theta)$ are the so-called action-angle coordinates, and reduce the transport operator $u\cdot\nabla$   to the much simpler
from $\frac{1}{T(h)}\de_\theta$, where $T(h)$, which is a $C^1$ function, is the period of the closed orbit $\{H(x_1,x_2)=h\}$.  The analogy with \eqref{eq:sheartransport} is then apparent, 
and estimate \eqref{eq:reg} is  derived from the explicit solution, obtained via the method of characteristics. Thanks to interpolation, the growth \eqref{eq:reg} is a lower bound on the
mix-norm of $\rho$, hence proving that $1/t$ is a lower bound on the mixing rate of $u$.

\subsection{Exponentially mixing flows}
Obtaining a faster mixing rate necessarily involves non-autonomous velocity fields. 
A widely used  exponential mixer, especially in numerical simulations, is due to Pierrehumbert \cite{pierrehumbert1994tracer}, and consists
of randomly alternating shear flows on $\T^2$. The beauty of this example is its simplicity: at discrete time steps $t_n$, it alternates the horizontal shear $(\sin (y-\omega_{1,n}),0)$
and the vertical shear  $(0,\sin (x-\omega_{2,n}))$.
Here, $\omega=\{\omega_{1,n},\omega_{2,n}\}$ is a sequence of independent uniformly distributed random variables 
so the phases are randomly shifted.
While widely believed to be exponentially mixing, the first proof of this fact appeared only recently in \cite{BCZG22}. 
\begin{theorem}\label{thm:BCZGmain}
There exists a random constant $D$ (with good bounds on its moments) and  $\gamma> 0$ such that we have~\eqref{e:expmix}, almost surely.
\end{theorem}

By taking a realization of the above velocity field, this result settles the question of the existence of a smooth exponential mixer on $\T^2$, although it does not produce a time-periodic velocity field.

The proof of Theorem \ref{thm:BCZGmain} relies on tools from random dynamical system theory and adopts a Lagrangian approach to the problem: this approach involves proving the positivity of the top Lyapunov exponent of the flow map
via Furstenberg's criterion and a Harris theorem.

A related example has been produced in \cites{HillSturmanEA22,ELM23}, constructed by alternating two piecewise linear shear flows. This example is fully deterministic and produces a time-periodic, Lipschitz velocity field. The important feature of this flow is that it generates a uniformly hyperbolic map on $\TT^2$.

In general, constructing exponentially mixing flows on $\T^2$ has proven to be quite a challenge, and only recently there have been tremendous developments. Besides the two works described above, that constitute the latest works in the field, we mention the deterministic constructions of \cites{AlbertiCrippaEA19,YaoZlatos17,ElgindiZlatos19}, and the beautiful work on velocity fields generated by stochastically forced Navier-Stokes equations of \cite{BBPSmix22}.

\section{Loss of regularity} \label{s:reg}

One of the effects of mixing is the creation of striation in the scalar field. Quantitatively, this effect corresponds to growth of derivatives of $\rho$, which can be seen from the interpolation inequality:
\begin{equation} \label{e:L2interpolation}
  \|\rho(t)\|_{L^2}^2 \leq \|\rho(t)\|_{\Dot{H}^{-s}} \|\rho(t)\|_{\Dot{H}^{s}}, \quad s>0. 
\end{equation}
In fact, recalling that the $L^2$ norm of $\rho$ is conserved by the flow of $u$, if the negative Sobolev norms of $\rho(t)$ decay to zero at some time $T\leq \infty$, at that same time the positive Sobolev norms must blow up. However, we note that growth of derivatives is a local phenomenon that can occur in the absence of mixing, which is a global phenomenon.

One can ask whether the growth of Sobolev norms can lead to loss of regularity for solutions of the transport equation \eqref{e:transport}, when the velocity field is not sufficiently smooth. 
The Cauchy-Lipschitz theory implies that, if $u$ is Lipschitz uniformly in time,  then the flow of $u$ is also Lipschitz continuous, although  its Lipschitz constant can grow exponentially fast in time. Hence, at least some regularity of the initial data $\rho_0$ is preserved in time. When the gradient of $u$ is not bounded, but it is still in some $L^p$ space with $p<\infty$, the direct estimates from the Cauchy-Lipschitz theory do not apply. Therefore, it is natural to investigate what, if any, regularity of the initial data $\rho_0$ is preserved under advection by $u$.

We present two examples to show that no Sobolev regularity is preserved in general: the first where the flow is mixing and all Sobolev regularity, including fractional regularity, is lost instantaneously;  
the second where the flow is not mixing and we are able to show at least that the $H^1$ (and any higher) norm blows up instantaneously. The second construction applies to (almost) all initial conditions in $H^1$, though the resulting velocity field $u$ still depends on the initial condition $\rho_0$. In both examples, the simple key idea is to utilize the linearity of the transport equation to construct a weak solutions by adding  infinitely-many suitably rescaled copies of a base flow and a base solution. The rescaling pushes energy to higher and higher frequencies or small scales, leading to an accelerated growth of the derivatives, which ultimately results in an instantaneous blow-up (see \cite{CrippaElgindiEA22} and references therein). 

The first result is the following: 

\begin{theorem} \label{t:ls1}
There exists a bounded velocity field $u$ such that $\nabla u(t) \in L^p(\RR^2))$, $1\leq p<\infty$, uniformly in time and a smooth, compactly supported  function $\rho_0\in C^\infty_c(\R^2)$, such that both $u$ and the unique bounded weak solution $\rho$ with initial data $\rho_0$ are compactly supported in space and smooth outside a point in $\RR^2$, but $\rho(t)$ does not belong to
$\dot{H}^s(\R^2)$ for any $s>0$ and $t>0$.
\end{theorem}


This result implies lack of continuity of the flow map in Sobolev spaces and can be shown to be a generic phenomenon in the sense of Baire's Category Theorem.

We sketch the proof of Theorem \ref{t:ls1}. Utilizing a suitable exponentially mixing flow on the torus, it is possible to construct a smooth, bounded, divergence-free vector field $u^0$ with $\nabla u^0(t) \in L^p(\RR^2)$, $1\leq p<\infty$, uniformly in time  and a smooth solution $\rho^0$ of the transport equation 
with velocity field $u^{(0)}$, both supported on the unit square $\mathcal{Q}_0$ in the plane, such that all positive Sobolev norms $\|\rho^0(t)\|_{\Dot{H}^s}$, $s>0$, grow exponentially fast in time.  For each $n\in \NN$, we define velocity fields $u^{(n)}$ and functions $\rho^{(n)}$ on squares $\mathcal{Q}_n$ of sidelength $\lambda_n$ by rescaling:
\begin{align} 
  u^{(n)}(t,x)= \frac{\lambda_n}{\tau_n}u^{(0)} \paren[\Big]{ \frac{t}{\tau_n},\frac{x}{\lambda_n}}, \nonumber \\
  \rho^{(n)}(t,x) = \gamma_n \rho^{(0)} \paren[\Big]{ \frac{t}{\tau_n},\frac{x}{\lambda_n}}, \label{e:lr1Rescaling}
\end{align}
 for some sequences $\lambda_n$, $\tau_n$, and $\gamma_n$ to be chosen,
up to some rigid motions, which do not change the norms and which we suppress for ease of notation. The squares $\mathcal{Q}_n$ can be taken pairwise disjoint. Then, by setting
\begin{equation} \label{e:Concatenate}
    u\defeq\sum_{n} u^{(n)}, \qquad \rho \defeq\sum_{n} \rho^{(n)},
\end{equation}
we have that $\rho$ is a weak solution of \eqref{e:transport} with velocity field $u$.
Lastly, we pick $\lambda_n$, $\tau_n$, and $\gamma_n$ in such a way that the squares $\mathcal{Q}_n$ converge to a point, the only point where $u$ and $\rho$ are not smooth, the norms  $\|u(t)\|_{\Dot{W}^{1,p}}$ and $\|\rho(0)\|_{\Dot{H}^s}$ are controlled, while the norm $\|\rho(t)\|_{\Dot{H}^s}=\infty$.

The second result is the following: 

\begin{theorem} \label{t:ls2}
 Given any non-constant function \mbox{$\rho_0 \in H^1(\RR^2)$}, there exists a bounded, compactly supported, divergence-free velocity field $u$, with $\nabla u(t) \in L^p(\RR^2)$ for any $1\leq p <\infty$ uniformly in time, and smooth outside a point in $\RR^2$, such that the unique weak solution $\rho(t)$ of~\eqref{e:transport} in $L^2(\RR^2)$ with initial data $\rho_0$ does not belong to  $\Dot{H}^1(\RR^2)$ (even locally) for any $t>0$.
\end{theorem}

In fact, a stronger statement is true. The velocity field $u$ is in all Sobolev spaces that are not embedded in the Lipschitz space (namely,  $u\in W^{r,p}$ for all $r<2/p+1$, $1\leq p<\infty$).

The main steps in the proof of Theorem \ref{t:ls2} are as follows. The first step follows by a direct calculation on functions on the torus $\TT^2$. Given any non-constant periodic function~$\bar \phi$, applying either a sine or cosine shear flow parallel to one of the coordinate axes must increase the~$\Dot{H}^1$ norm of~$\bar\phi$ by a constant factor at time $t=1$. Hence, the $\Dot{H}^1$-norm of $\bar\phi$ grows exponentially in time. By unfolding the action of the shear flows on the torus, this observation can be adapted to showing exponential growth of the $\Dot{H}^1$-norm of functions supported on a square (let's say again the unit square) in $\RR^2$ by a combination of sine and cosines shear flows. Next, rescaling the flow alone in a manner similar to \eqref{e:lr1Rescaling} gives a sequence of well-separated squares, shrinking to a point, such that the rescaled flow grows the $\Dot{H}^1$-norm of functions supported on each square by a larger and larger factor at time 1. The final step consists in choosing the location of the squares and the rescaling factor in such a way that the $\Dot{H}^1$ norm of the solution diverges at any positive time, but the velocity field remains sufficiently regular. This last step involves a certain covering lemma to ensure that the gradient of the initial data is sufficient large in $L^2$ averaged sense on all the squares.

 In view of the results above, one can ask if any regularity of $\rho_0$, measured by a norm that is not comparable with the Sobolev norm $\Dot{H}^s$, is preserved under the advection by $u$. It can be shown that essentially only the ``logarithm" of a derivative
 is preserved \cite{BrueNguyen21}. In this sense, the loss of regularity in Theorem \ref{t:ls1}  can be viewed as optimal.

\section{Enhanced dissipation} \label{s:ed}

We now turn our attention to problems where both diffusion and convection are present, and study the combined effect of both. 
A prototypical example is the evolution of the concentration of a solute in an ambient fluid (e.g.~cream in coffee).
The evolution of the (normalized) concentration is modelled by the advection diffusion equation
\begin{equation}\label{e:AdvectionDiffusion}
  \partial_t \rho^\kappa + u \cdot \grad \rho^\kappa - \kappa \lap \rho^\kappa = 0\,.
\end{equation}
Here~$\rho^\kappa$ denotes the deviation of the concentration of the solute from its spatial average, the quantity $\kappa > 0$ is the molecular diffusivity, and $u$ denotes the velocity field of the ambient fluid.
As in the previous sections, we will impose periodic boundary conditions on~\eqref{e:AdvectionDiffusion} and restrict our attention to the incompressible setting where~$u$ is divergence free.
In this case the spatial average is still preserved by~\eqref{e:AdvectionDiffusion}, so we may, without loss of generality, uphold our convention that $\rho^\kappa$ is spatially mean zero.
%

Multiplying~\eqref{e:AdvectionDiffusion} by~$\rho$, integrating in space, using incompressibility and  Poincar\'e's inequality shows
\begin{equation}\label{e:Poincare}
  \norm{\rho^\kappa(t)}_{L^2} \leq \e^{ - 4 \pi^2 \kappa  t} \norm{\rho_0}_{L^2}
  \,.
\end{equation}
This estimate, however, is completely blind to the effect of the fluid advection, and in practice one expects~$\norm{\rho^\kappa}_{L^2}$ to decay much faster than the rate predicted by~\eqref{e:Poincare}.
The reason for this is a phenomenon most people have likely observed themselves when stirring cream into coffee: the fluid flow initially spreads the cream into fine filaments; diffusion acts faster on fine filaments, and so these uniformize very quickly.

\subsection{Quantifying dissipation enhancement.}
One way to mathematically quantify and study this phenomenon is through the \emph{dissipation time}, denoted by~$\tdis$.
Explicitly define $\tdis = \tdis(u, \kappa)$ to be the smallest time $t \geq 0$ so that
\begin{equation}
  \norm{\rho^\kappa(s + t)}_{L^2} \leq \frac{1}{2} \norm{\rho^\kappa(s)}_{L^2}\,,
\end{equation}
for every time $s \geq 0$ and mean-zero initial data $\rho^\kappa(s) \in L^2$.
Clearly~\eqref{e:Poincare} shows~$\tdis \leq 1 / (4 \pi^2 \kappa)$, and one can precisely define \emph{enhanced dissipation} as situations where~$\tdis \ll 1 / (4 \pi^2 \kappa)$.
We now list several situations where enhanced dissipation is exhibited.

\emph{Shear flows.}
If~$u$ is a shear flow with a~$C^2$ profile that has non-degenerate critical points, then classical work of Kelvin shows~$\tdis \leq C / \sqrt{\kappa}$.
A matching lower bound was also proved by Coti Zelati and Drivas.

\emph{Cellular flows.}
A cellular flow models the movement of a 2D fluid in the presence of a strong array of opposing vortices.
The simplest example is given by
\begin{equation}
  u = \begin{pmatrix}
    -\partial_{x_2} H
    \\
    \partial_{x_1} H
  \end{pmatrix}\,,
  \quad
  H = A \epsilon \sin\paren[\Big]{\frac{2 \pi x_1}{\epsilon}}
    \sin\paren[\Big]{\frac{2 \pi x_2}{\epsilon}}\,,
\end{equation}
where~$A \gg 1$ is the flow amplitude and $\epsilon \ll 1$ is the cell size.
Standard homogenization results show that as $\epsilon \to 0$, the operator $-u \cdot \grad + \kappa \lap$ 
behaves like $D_\mathrm{eff} \lap$,
where
\begin{equation}
  D_\mathrm{eff} \approx C \sqrt{\frac{\kappa A}{\epsilon} }\,,
\end{equation}
is the effective diffusivity.
As a result one would expect
\begin{equation}
  \tdis
    = O\paren[\Big]{ \frac{1}{D_\mathrm{eff}} }
    = O\paren[\Big]{ \sqrt{\frac{\epsilon}{\kappa A} } }
\end{equation}
as~$\kappa, \epsilon \to 0$, $A \to \infty$.
This is indeed the case (provided $\kappa / \epsilon \ll A \ll \kappa / \epsilon^3$), and was proved recently by Iyer and Zhou.
A matching lower bound for~$\tdis$ was recently proved in~\cite{BCZM22}.

\emph{Relaxation enhancing flows.}
Seminal work of Constantin \etal~\cite{ConstantinKiselevEA08} shows that for time independent velocity fields, $\tdis = o(1/\kappa)$ if and only if the operator~$u \cdot \grad$ has no eigenfunctions in~$H^1$.
Such flows are called \emph{relaxation enhancing}.
It is known that \emph{weakly mixing} flows are relaxation enhancing, but relaxation enhancing flows need not be weakly mixing.

\emph{Mixing flows.}
Thus far, the examples provided only reduce the dissipation time to an algebraic power~$1/\kappa^\alpha$ for some $\alpha < 1$.
Using \emph{exponentially mixing} flows, it is possible to reduce the dissipation time further to~$\abs{\ln \kappa}^2$
(see~\cites{FengIyer19,CotiZelatiDelgadinoEA20}).
In fact, we recall that a velocity field~$u$ is \emph{exponentially mixing} if, for any~$\rho_0 \in \dot H^1$, the mix norm of solutions to~\eqref{e:transport} decays exponentially as in~\eqref{e:expmix}.
Informally, then any solution to~\eqref{e:transport} with initial data that is localized to a ball of radius~$\epsilon$ becomes essentially uniformly spread in time $O(\abs{\ln \epsilon})$.
If~$u$ is exponentially mixing, one can show that
\begin{equation}\label{e:lka2}
  \tdis \leq C \abs{\ln \kappa}^2\,.
\end{equation}
An elementary heuristic argument, however, suggests we should have the stronger bound
\begin{equation}\label{e:lka}
  \tdis \leq C \abs{\ln \kappa}\,.
\end{equation}
Indeed, if the solute is initially concentrated at one point~$x$, then after time $O(1)$ it will spread to a ball of radius~$O(\sqrt{\kappa})$.
Now, since~$u$ is exponentially mixing, it will get spread  almost uniformly on the entire domain in time~$O(\abs{\ln \kappa})$, if the effect of diffusion is negligible.
Unfortunately, the effects of diffusion may not be negligible on the time scales of order~$O(\abs{\ln \kappa})$, and so this argument cannot be easily made rigorous.

Even though proving~\eqref{e:lka} for general exponentially mixing flows is still an open question, there are several examples of flows for which~\eqref{e:lka} is known: for instance, when~$u$ is the velocity field from the stochastically forced Navier--Stokes equations~\cite{BedrossianBlumenthalEA21}, or when $u$ consists of  alternating horizontal/vertical shears with a tent profile and a
sufficiently large amplitude~\cite{ELM23}.
It is also known that $\tdis$ cannot be smaller than $O(\abs{\ln \kappa})$ for velocity fields that are Lipschitz in space uniformly in time.

\subsection{Blow up suppression}

One application of enhanced dissipation is to control certain nonlinear phenomena.
For concreteness and simplicity, we focus our attention on a simplified version of the Keller--Segel system of equations, which is used to model the evolution of the population density of micro-organisms when chemotactic effects are present.
We again impose periodic boundary conditions and study this system on $2$  dimensional torus.
If~$n = n(t, x)\geq 0$ represents the bacterial population density, ~$c\geq 0$ represents the concentration of a chemoattractant  produced by the bacteria,
and $\chi > 0$ is a sensitivity parameter, then a simplified version of the Keller--Segel model is the following system:
\begin{gather}
  \label{e:ks1}
  \partial_t n - \lap n = -\dv \paren[\big]{ n \chi \grad c }\,,
  \\
  \label{e:ks2}
  - \lap c = n - \bar n\,,
  \\
  \label{e:ks3}
  \bar n = \int_{\T^2} n \, \dd x\,.
\end{gather}
This model stipulates that bacterial diffusion is biased in the direction of the gradient of the concentration of a chemoattractant that is emitted by the bacteria themselves, and that the chemoattractant diffuses much faster than the bacteria do.

From the equation we see that there is a competition between two effects: The diffusive term~$\lap n$ drives bacteria away from regions of high population, and the chemotactic term $\dv (n \chi \grad c)$ drives the bacteria towards it.
If the chemotactic effects dominate, they will lead to a population explosion.
This has been well studied and it is now known that the diffusive effects dominate (and there is no population explosion) if and only if the total initial  population is below a certain threshold.

One natural question is to study the effect of movement of the ambient fluid on this system.
If the ambient fluid has velocity field~$u$, 
equation~\eqref{e:ks1} becomes
\begin{equation}\label{e:ks1drift}
  \tag{\ref{e:ks1}$'$}
  \partial_t n + u \cdot \grad n - \lap n = -\dv \paren[\big]{ n \chi \grad c }\,.
\end{equation}
Intuitively, we expect that regions of high concentration of bacteria can be dispersed by vigorous stirring.
This result can be established rigorously.

\begin{theorem}\label{t:noblowup}
  There exists~$t_* = t_*( \norm{n_0}_{L^2})$ such that if
  \begin{equation}
    \tdis(u, 1) < t_*\,,
  \end{equation}
  then there is no population explosion in~\eqref{e:ks1drift}.
\end{theorem}

The main idea being the proof can be explained in an elementary fashion and we now provide a quick sketch of the proof of Theorem~\ref{t:noblowup}.

Next, we rewrite~\eqref{e:ks1drift} as
  \begin{equation}\label{e:KSTheta}
    \partial_t \theta + u \cdot \grad \theta - \lap \theta = \mathcal N(\theta)\,,
  \end{equation}
where $\theta = n - \bar n$ and
  \begin{equation}
    \mathcal N(\theta) = - \dv( (\theta + \bar n) \chi \grad c ) \,.
  \end{equation}
  Multiplying~\eqref{e:KSTheta} by~$\theta$ and integrating in space show that 
  \begin{equation}
    \frac{1}{2} \partial_t \norm{\theta(t)}_{L^2}^2 + \norm{\grad \theta}_{L^2}^2 \leq \int_{\T^d} \theta \mathcal N(\theta) \, \dd x \,.
  \end{equation}
  Now using standard energy estimates one can show that when the diffusive term $\norm{\grad \theta}_{L^2}^2$ is large, then $\partial_t \norm{\theta}_{L^2} \leq 0$.
  More precisely, one can show there exists a constant $C_1 = C_1( \chi, d, \norm{\theta_0}_{L^2} )$ such that if
  \begin{equation}\label{e:H1Large}
    \frac{1}{t_*} \int_0^{t_*} \norm{\grad \theta}_{L^2}^2 \, \dd t \geq C_1\,,
  \end{equation}
  then we must also have~$\norm{\theta(t_*)}_{L^2} \leq \norm{\theta_0}_{L^2}$.

  Suppose now~\eqref{e:H1Large} does not hold.
  In this case we use Duhamel's formula to write
  \begin{equation}
    \theta(t_*) = \mathcal S_{0, t_*} \theta_0  + \int_0^{t_*} S_{s, t_*} \mathcal N(\theta(s)) \, ds\,,
  \end{equation}
  where $\mathcal S_{s, t}$ is the solution operator to~\eqref{e:AdvectionDiffusion}.
  This implies
  \begin{equation}
    \norm{\theta(t_*)}_{L^2}^2 \leq \norm{\mathcal S_{0, t_*} \theta_0}_{L^2}^2  + \int_0^{t_*} \norm{ \mathcal N(\theta(s))}_{L^2} \, ds\,.
  \end{equation}

  Notice that, since~$\tdis \leq t_*$, the first term on the right is at most $\norm{\theta_0}_{L^2}/ 2$.
  For the second term we use standard energy estimates to control $\norm{\mathcal N(\theta)}_{L^2}$ by $\norm{\grad \theta}_{L^2}^2$ and~$\norm{\theta}_{L^2}$.
  Combined with the assumption that~\eqref{e:H1Large} does not hold, we obtain an inequality  of the form
  \begin{equation}
    \norm{\theta(t_*)}_{L^2}^2 \leq \paren[\Big]{\frac{1}{2} + t_* F(C_1) } \norm{\theta_0}_{L^2}^2 \,,
  \end{equation}
  for some explicit function~$F$ that arises from the bound on the nonlinearity.
  If~$t_*  \leq 1 / (2 F(C_1))$ then the right hand side is at most~$\norm{\theta_0}_{L^2}$.
  Iterating this step, one immediately sees that $\sup_{t < \infty} \norm{\theta(t)}_{L^2} < \infty$.
  This bound is enough to show that there is no population explosion in~\eqref{e:ks1drift}, concluding the proof.
  A more complete version of this proof can be found in~\cite{IyerXuEA21}.
  It has also be extended to fourth order equations, such as the Kuramoto-Sivashinsky equation, a model of flame-front propagation, by Feng and Mazzucato.

\section{Anomalous dissipation} \label{s:ad}

In the previous section we saw several examples of \emph{enhanced dissipation}, where the solution to~\eqref{e:AdvectionDiffusion} loses a constant fraction of its $L^2$ energy in time scales much smaller than the dissipative time scale~$1/\kappa$.
The examples outlined exhibited the energy loss on time scales~$1/\kappa^\alpha$, $\abs{\ln \kappa}^2$ or~$\abs{\ln \kappa}$, which diverge to infinity in the vanishing diffusivity limit~$\kappa \to 0$.
A natural question to ask is whether there are situations where solutions to~\eqref{e:AdvectionDiffusion} lose a constant fraction of their $L^2$ energy on a time scale that is~$O(1)$ as~$\kappa \to 0$.
This phenomenon is called~\emph{anomalous dissipation}.
More precisely, anomalous dissipation is the existence of solutions of~\eqref{e:AdvectionDiffusion} with progressively smaller diffusivity $\kappa$ that converge to a \emph{dissipative} solution of the transport equation in the limit $\kappa\to0$.
That is, as~$\kappa \to 0$, we can find solutions $\rho^\kappa$ to~\eqref{e:AdvectionDiffusion},  which converge (possibly along a subsequence) to a solution $\rho^0$ of~\eqref{e:transport},  where
$$
\| \rho^0(T) \|_{L^2}^2
< \| \rho^0(0) \|_{L^2}^2 \,,
$$
for some time~$T < \infty$.

The dissipation of the $L^2$ energy of the solution $\rho^\kappa$ to~\eqref{e:AdvectionDiffusion} for time $t \in [0,T]$ is encoded in the energy estimate 
\begin{multline}\label{e:balance}
\| \rho^\kappa(t) \|_{L^2}^2
- \| \rho^\kappa(0) \|_{L^2}^2
= - 2 \kappa \int_0^t \| \nabla \rho^\kappa (s) \|_{L^2}^2 \, ds \\
\qquad \text{ for any $t \in [0,T]$.}
\end{multline}
In the limit case $\kappa=0$,~\eqref{e:balance} expresses (formally) the conservation of the $L^2$ norm for solutions of the transport equation~\eqref{e:transport}. The velocity field does not appear explicitly in~\eqref{e:balance}. However, the action of a mixing velocity field results in filamentation of the scalar and consequently in the creation of large gradients, thus conceivably allowing for scenarios in which the right hand side of~\eqref{e:balance} remains bounded away from zero even in the limit $\kappa \to 0$ of vanishing diffusivity.

This phenomenon is the analogue in the linear case of the so-called $0$-th law of turbulence of the Onsager-Kolmogorov theory of turbulence for the Euler/Navier-Stokes equations. The $0$-th law predicts uniform-in-viscosity dissipation of the kinetic energy, due to the nonlinear transfer of energy to high frequencies and to the corresponding enhanced effect of the diffusion.

In order to identify the critical regularity for anomalous dissipation in solutions to~\eqref{e:AdvectionDiffusion} as $\kappa$ vanishes, heuristically at least,  we can formally rewrite the contribution of the advection term in the energy estimate as 
$$
\int_{\T^2}  ( u \cdot \nabla \rho ) \rho \, \dd x 
\sim \int_{\T^2}  \nabla^\alpha u \, \big( \nabla^{\frac{1-\alpha}{2}} \rho \big)^2 \, \dd x \,,
$$
for any $0\leq\alpha\leq 1$, where $\rho$ denotes both $\rho^\kappa$ and $\rho^0$. The fractional derivatives $\nabla^\alpha$ can be estimated via norms in H\"older's spaces $C^\alpha$.
Criticality is therefore expressed by the so-called Yaglom's relation: focusing for simplicity on regularity in space only, for $u \in C^\alpha$ and $\rho \in C^\beta$, the combined H\"older's regularity of the velocity field and the solution is:
\begin{itemize}
    \item subcritical, if $\alpha + 2 \beta > 1$,
    \item critical, if $\alpha + 2 \beta =1$,
    \item supercritical, if $\alpha + 2 \beta <1 $.
\end{itemize}
This regularity is the analogue for linear advection-diffusion equations of the critical $\frac{1}{3}$-H\"older's regularity in the case of anomalous energy dissipation for solutions of the Navier-Stokes equations as viscosity vanishes, according to the Onsager-Kolmogorov theory of turbulence. This analogy can be formally seen by replacing the fluid  velocity in the Navier-Stokes and Euler equations  with $\rho$.

The Obukhov-Corrsin theory of scalar turbulence (1949-1951) predicts that:
\begin{itemize}
    \item in the subcritical regime, for a given $u \in C^\alpha$ there exists a unique solution $\rho^0\in C^\beta$ of~\eqref{e:transport} and such a solution conserves the $L^2$ norm,
    \item in the supercritical regime, there exist velocity fields $u\in C^\alpha$ such that nonuniqueness and dissipation of the $L^2$ norm are possible for solutions $\rho^0\in C^\beta$ of~\eqref{e:transport}; moreover, anomalous dissipation is possible, in the sense that
    \begin{equation}\label{e:lasttime}
    \limsup_{\kappa\to 0} \; \kappa \int_0^T \| \nabla \rho^\kappa(s) \|^2_{L^2} \, \dd s > 0 
    \end{equation}
    for solutions $\rho^\kappa$ of~\eqref{e:AdvectionDiffusion} uniformly bounded in~$C^\beta$.
    \end{itemize}
The statement in the subcritical case can be proven using a commutator argument similar to that of Constantin, E, and Titi for the Euler equations, see for instance~\cite{DrivasElgindiEA22}*{Theorem~4}. Notice that the uniqueness statement strongly relies on the linearity of the equation. Addressing the supercritical case is more challenging and has been done only very recently from a rigorous mathematical perspective. We briefly discuss these recent results in the next Subsection. 

\subsection{Anomalous dissipation for boun\-ded solutions}

The endpoint case $\alpha<1$ and $\beta=0$ has been addressed in~\cite{DrivasElgindiEA22}. In this paper, for any $\alpha <1$, the authors provide an example of a bounded velocity field, which belongs to $L^1([0,T];C^\alpha(\mathbb{T}^2))$ and is smooth except at the singular time $t=T$, that exhibits anomalous dissipation for all initial data sufficiently close to a (nontrivial) harmonics. They are also able to construct velocity fields that exhibit anomalous dissipation for any (regular enough) initial datum, although the velocity field depends on the chosen initial datum. The strategy for both examples is to construct a velocity field which develops smaller and smaller scales when the time approaches the singular time $t=T$. This construction can be interpreted as mimicking the development in time of a turbulent cascade. However, it also causes the anomalous dissipation to be concentrated at the singular time $t=T$, in the sense that for any $\varepsilon>0$ it holds
    $$
    \lim_{\kappa\to 0} \; \kappa \int_0^{T-\varepsilon} \| \nabla \rho^\kappa(s) \|^2_{L^2} \, \dd s = 0 \,.
    $$
Several criteria that imply anomalous dissipation are given in~\cite{DrivasElgindiEA22}. The criterion that provides the most intuition on the mechanism for anomalous dissipation (even though such a criterion is not the one effectively exploited in the proof in~\cite{DrivasElgindiEA22}) establishes a link to mixing in solutions of the transport equation~\eqref{e:transport} and asserts that, if the solution $\rho$ of~\eqref{e:transport} satisfies
\begin{multline}\label{e:inverse}
\int_0^T \| \nabla \rho(s) \|^2_{L^2} \, \dd s = \infty
\\  \text{ and } \quad  
\| \rho(t) \|_{\dot{H}^{-1}} \| \rho(t) \|_{\dot{H}^{1}} 
\leq C \| \rho(t) \|^2_{L^2}\,,
\end{multline}
then anomalous dissipation holds. 
In view of the interpolation inequality \eqref{e:L2interpolation} with $s=1$, the second condition in~\eqref{e:inverse} in particular implies that, in the absence of diffusion, the velocity field mixes essentially at the optimal rate.
However, it is not easy to produce velocity fields with such strong mixing properties as in~\eqref{e:inverse}, and this is the reason why the proof of~\cite{DrivasElgindiEA22} needs to rely on weaker criteria. In fact,  the velocity field is a self-similar version of the alternating shear flows example by Pierrehumbert~\cite{pierrehumbert1994tracer}.  Although it enjoys  weaker mixing properties, the velocity field constructed in~\cite{DrivasElgindiEA22} exhibits  anomalous dissipation due to the following heuristic reason: mixing requires all energy to be sent to high frequencies, while anomalous dissipation just requires a given fraction of the energy to be sent to high frequencies. 

As a corollary, the construction in~\cite{DrivasElgindiEA22} provides a new example of nonuniqueness for the transport equation~\eqref{e:transport} with a velocity field in $L^1([0,T];C^\alpha(\mathbb{T}^2))$ (for $\alpha<1$), but outside the DiPerna-Lions class. In order to see this fact, we extend the velocity field for time $t \in [T,2T]$ by reflecting it oddly across $t=T$, that is, setting $u(t) = - u(2T-t)$. We can then see that the following are two distinct weak solutions: 
\begin{itemize}
    \item the vanishing-diffusivity solution $\rho^{\rm{vd}}(t)$ for time $t\in[0,2T]$, and
    \item the solution given by $\rho^{\rm{refl}} (t) = \rho^{\rm{vd}}(t)$ for time $t\in[0,T]$, and by the odd reflection $\rho^{\rm{refl}} (t) = \rho^{\rm{vd}}(2T-t)$ for time $t \in [T,2T]$.
\end{itemize}
Indeed, they are distinct for $t \in [T,2T]$, since $\rho^{\rm{vd}}$ dissipates the~$L^2$ norm while $\rho^{\rm{refl}}$ conserves it. It must be noted that, based on the approach in~\cite{DrivasElgindiEA22}, it is unclear whether $\rho^{\rm{refl}}$ can be constructed in the limit of vanishing diffusivity along a suitable subsequence $\kappa_n \to 0$.

\subsection{Anomalous dissipation and lack of selection}

The possibility of having two distinct solutions, both arising in the limit of vanishing diffusivity, puts into question the validity of the zero-diffusivity limit as a selection principle for weak solutions of the transport equation outside of the DiPerna--Lions theory. Solutions arising as zero-diffusivity limit may be considered ``more physical'' than general weak solutions and one may wonder whether uniqueness could be restored in the sense of the possibility of a selection mechanism among the many weak solutions. This was the leading question behind the analysis in~\cite{CCS}, in which it was shown that neither vanishing diffusivity nor regularization of the velocity field provide such a selection mechanism. 

In terms of anomalous dissipation, in \cite{CCS} the authors construct, for any $\alpha$ and $\beta$ in the supercritical regime $\alpha+2\beta<1$,  a velocity field $u \in L^\infty([0,T];C^\alpha(\mathbb{T}^2))$ and a smooth initial datum $\rho_0$ such that the solutions $\rho^\kappa$ to the advection-diffusion equation~\eqref{e:AdvectionDiffusion} with initial data $\rho_0$ are uniformly bounded in $L^2([0,T];C^\beta(\mathbb{T}^2))$ and exhibit anomalous dissipation (more general exponents for the integrability in time  can be considered). 

The basic mechanism is based on the same slice-and-dice strategy that leads to the non-uniqueness results of Depauw and Bressan for the transport equation and is intrinsically related to mixing. In such an example, the solution takes opposite constant values on alternating tiles of a checkerboard, the size of which gets refined as time increases, resulting in perfect mixing (i.e.,~weak convergence to zero, the average of the solution) at the critical time $t=T$. Compared to the previous literature on mixing, a few important twists are necessary for the analysis in \cite{CCS}. The refinement of the size of the checkerboards does not follow the classical dyadic rescaling, but rather obeys a superexponential scaling law, which allows (upon suitable choices of the many parameters in the construction) to achieve optimal regularity and to separate the relevant scales at each step of the evolution. The transition from a checkerboard to the next one is realized by shear flows, which are also concentrated at suitable spatial scales.
The approach of \cite{CCS} is fairly explicit and relies on the fact that solutions of the advection-diffusion equation~\eqref{e:AdvectionDiffusion} have a stochastic Lagrangian representation
via the Feynman-Kac formula
$$
\rho^\kappa(t,x) = \mathbb{E} \big[ \rho_0 \big((X^\kappa)^{-1}(t,x) \big) \big] \,.
$$
Here~$X^\kappa$ satisfies the stochastic differential equation
$$
d X^\kappa (t,x) = u(t,X^\kappa(t,x)) \, \dd t + \sqrt{2\kappa} \, dW\,,
$$
where $W$ is a Brownian motion. 

The principle behind the lack of selection in the limit of zero diffusivity  can be best understood by first considering the related question of whether regularizing the velocity via convolution with a suitable smoothing kernel can act as a selection principle for the solution, as the regularization parameter vanishes. That is, we pose the question whether limit points of solutions $\rho^\varepsilon$ of
$$
\de_t\rho^\varepsilon + (u\ast \eta_\varepsilon) \cdot\nabla \rho^\varepsilon=0
$$
are unique, where $\eta_\varepsilon$ is a standard mollifier. To this extent, consider the slide-and-dice construction for $t \in [0,T]$ sketched above, and reflect it oddly across $t=T$, therefore ``reconstructing large scales'' for $t\in[T,2T]$. At each step in the reconstruction of the large scales, add a new ``move'' which has the effect of changing the parity of the checkerboard, that is, it swaps the black and the white tiles. When considering the convolution of the velocity field with $\eta_\varepsilon$, all scales below $\varepsilon$ are ``filtered'' and therefore the solution does not get fully mixed at time $t=T$, but rather stays at a finite scale when it crosses the singular time. After the singular time, large scales are reconstructed, but, as $\varepsilon \to 0$, the solution will converge  either to an even or odd checkerboard for $t \sim 2T$,  depending on the ``parity of~$\varepsilon$''. This construction provides two subsequences that converge to distinct limit solutions. 


We next focus on the case when diffusion is present. Consider again the evolution of the checkerboards for $t \in [0,T]$. Now, add at each step a time interval of suitable length on which the velocity vanishes and, therefore, the solution obeys the heat equation. Separation of scales (due to the choice of a superexponential sequence) provides the existence of a critical time $t_{\rm crit}(\kappa) \to T$ such that for $0<t< t_{\rm crit}(\kappa)$ diffusion is a perturbation, while it is the dominant effect for times $t \sim t_{\rm crit}(\kappa)$ due to the intervals where velocity vanishes. Although acting only for a short time, diffusion is enhanced by the high frequencies in the solution, eventually reaching a balance which leads to dissipate a fixed fraction, independent of $\kappa>0$,  of the $L^2$ norm of the solution. Hence, anomalous dissipation occurs. 

In the construction sketched above, the possibility of anomalous dissipation relies on a specific choice of the subsequence $\kappa \to 0$ that depends on all other parameters. In fact, another result of \cite{CCS} is the possibility of choosing another (distinct) subsequence $\kappa \to 0$ with the following property.
Exploiting the isotropy of the Brownian motion, the corresponding subsequence of solutions $\rho^\kappa$ converges to a solution of the transport equation which conserves the~$L^2$ norm. This example shows that the limit of zero diffusion cannot be used as a selection principle for weak solutions of the transport equation outside of the DiPerna-Lions class. 

Building on the results of \cite{CCS}, in collaboration with Bru\`e and De Lellis the authors have obtained analogous results for the forced Navier-Stokes and Euler equations with full Onsager-supercritical regularity, i.e., for velocity fields in H\"older spaces $C^\alpha$ with any $\alpha<1/3$.

\subsection{Anomalous dissipation via fractal homogeneization}

Both in~\cite{DrivasElgindiEA22} and in~\cite{CCS},  anomalous dissipation only occurs at the singular time $t=T$ (see~\eqref{e:lasttime}), due to the nature of the construction based on mixing and on the development in time of small scales. Such a situation is somewhat inconsistent with the theory of homogeneous isotropic turbulence, which postulates (statistical) stationarity and, therefore, the fact that there should be no ``preferred'' time in turbulent phenomena: anomalous dissipation should happen continuously in time, and for any randomly chosen times $t_1$ and $t_2$ the corresponding values of the velocity field $u(t_1,\cdot)$ and $u(t_2,\cdot)$ should be macroscopically indistinguishable. 

In~\cite{AV}, the authors address this issue by relying on an approach based on homogenization theory. For any $\alpha < \frac{1}{3}$, they construct a time-periodic velocity field on the torus that belongs to $C^\alpha$ uniformly in time and exhibits continuous-in-time anomalous dissipation for bounded solutions with arbitrary $H^1$ initial data.

In contrast to the examples in \cite{CCS} and~\cite{DrivasElgindiEA22}, for the velocity field constructed in~\cite{AV} all scales are active at all times.
Homogenization theory allows to understand and quantify the enhancement of diffusivity due to the creation of small scales in the solution at all times. This result goes under the name of renormalization of effective diffusivities: each homogenization step along the cascade of scales enhances the effective diffusivity, which after an iteration over all scales remains of order one even as $\kappa \to 0$ and, therefore, gives anomalous dissipation.
Homogenization is well understood for a finite number of scales, but the authors of~\cite{AV} need to deal with infinitely many scales at once.

\bibliography{gautam,anna}

\begin{bibdiv}
\begin{biblist}

\bib{AlbertiCrippaEA19}{article}{
      author={Alberti, Giovanni},
      author={Crippa, Gianluca},
      author={Mazzucato, Anna~L.},
       title={Exponential self-similar mixing by incompressible flows},
        date={2019},
        ISSN={0894-0347},
     journal={J. Amer. Math. Soc.},
      volume={32},
      number={2},
       pages={445\ndash 490},
      review={\MR{3904158}},
}

\bib{Ambrosio04}{article}{
      author={Ambrosio, Luigi},
       title={Transport equation and {C}auchy problem for {$BV$} vector
  fields},
        date={2004},
        ISSN={0020-9910},
     journal={Invent. Math.},
      volume={158},
      number={2},
       pages={227\ndash 260},
      review={\MR{2096794}},
}

\bib{AV}{misc}{
      author={Armstrong, Scott},
      author={Vicol, Vlad},
       title={Anomalous diffusion by fractal homogenization},
        date={2023},
}

\bib{BedrossianBlumenthalEA21}{article}{
      author={Bedrossian, Jacob},
      author={Blumenthal, Alex},
      author={Punshon-Smith, Sam},
       title={Almost-sure enhanced dissipation and uniform-in-diffusivity
  exponential mixing for advection-diffusion by stochastic {N}avier-{S}tokes},
        date={2021},
        ISSN={0178-8051},
     journal={Probab. Theory Related Fields},
      volume={179},
      number={3-4},
       pages={777\ndash 834},
      review={\MR{4242626}},
}

\bib{BBPSmix22}{article}{
      author={Bedrossian, Jacob},
      author={Blumenthal, Alex},
      author={Punshon-Smith, Samuel},
       title={Almost-sure exponential mixing of passive scalars by the
  stochastic {N}avier-{S}tokes equations},
        date={2022},
        ISSN={0091-1798,2168-894X},
     journal={Ann. Probab.},
      volume={50},
      number={1},
       pages={241\ndash 303},
         url={https://doi.org/10.1214/21-aop1533},
      review={\MR{4385127}},
}

\bib{BCZM22}{article}{
      author={{Bru{\`e}}, Elia},
      author={{Coti Zelati}, Michele},
      author={{Marconi}, Elio},
       title={{Enhanced dissipation for two-dimensional Hamiltonian flows}},
        date={2022-11},
     journal={arXiv e-prints},
      eprint={2211.14057},
}

\bib{BCZG22}{article}{
      author={Blumenthal, Alex},
      author={Coti~Zelati, Michele},
      author={Gvalani, Rishabh~S.},
       title={Exponential mixing for random dynamical systems and an example of
  {P}ierrehumbert},
        date={2023},
        ISSN={0091-1798,2168-894X},
     journal={Ann. Probab.},
      volume={51},
      number={4},
       pages={1559\ndash 1601},
         url={https://doi.org/10.1214/23-aop1627},
      review={\MR{4597327}},
}

\bib{BrueNguyen21}{article}{
      author={Bru\`e, Elia},
      author={Nguyen, Quoc-Hung},
       title={Sharp regularity estimates for solutions of the continuity
  equation drifted by {S}obolev vector fields},
        date={2021},
        ISSN={2157-5045,1948-206X},
     journal={Anal. PDE},
      volume={14},
      number={8},
       pages={2539\ndash 2559},
         url={https://doi.org/10.2140/apde.2021.14.2539},
      review={\MR{4377866}},
}

\bib{Bressan03}{article}{
      author={Bressan, Alberto},
       title={A lemma and a conjecture on the cost of rearrangements},
        date={2003},
        ISSN={0041-8994},
     journal={Rend. Sem. Mat. Univ. Padova},
      volume={110},
       pages={97\ndash 102},
      review={\MR{2033002 (2005a:35187)}},
}

\bib{CCS}{misc}{
      author={Colombo, Maria},
      author={Crippa, Gianluca},
      author={Sorella, Massimo},
       title={Anomalous dissipation and lack of selection in the
  {O}bukhov-{C}orrsin theory of scalar turbulence},
        date={2022},
}

\bib{CrippaDeLellis08}{article}{
      author={Crippa, Gianluca},
      author={De~Lellis, Camillo},
       title={Estimates and regularity results for the {D}i{P}erna-{L}ions
  flow},
        date={2008},
     journal={J. Reine Angew. Math},
      volume={616},
       pages={15\ndash 46},
}

\bib{CrippaElgindiEA22}{article}{
      author={Crippa, Gianluca},
      author={Elgindi, Tarek},
      author={Iyer, Gautam},
      author={Mazzucato, Anna~L.},
       title={Growth of {S}obolev norms and loss of regularity in transport
  equations},
        date={2022},
        ISSN={1364-503X,1471-2962},
     journal={Philos. Trans. Roy. Soc. A},
      volume={380},
      number={2225},
       pages={Paper No. 24, 12},
         url={https://doi.org/10.1098/rsta.2021.0024},
      review={\MR{4430388}},
}

\bib{ConstantinKiselevEA08}{article}{
      author={Constantin, P.},
      author={Kiselev, A.},
      author={Ryzhik, L.},
      author={Zlato{\v{s}}, A.},
       title={Diffusion and mixing in fluid flow},
        date={2008},
        ISSN={0003-486X},
     journal={Ann. of Math. (2)},
      volume={168},
      number={2},
       pages={643\ndash 674},
      review={\MR{2434887 (2009e:58045)}},
}

\bib{CotiZelatiDelgadinoEA20}{article}{
      author={Coti~Zelati, Michele},
      author={Delgadino, Matias~G.},
      author={Elgindi, Tarek~M.},
       title={On the relation between enhanced dissipation timescales and
  mixing rates},
        date={2020},
        ISSN={0010-3640},
     journal={Comm. Pure Appl. Math.},
      volume={73},
      number={6},
       pages={1205\ndash 1244},
      review={\MR{4156602}},
}

\bib{DrivasElgindiEA22}{article}{
      author={Drivas, Theodore~D.},
      author={Elgindi, Tarek~M.},
      author={Iyer, Gautam},
      author={Jeong, In-Jee},
       title={Anomalous dissipation in passive scalar transport},
        date={2022},
        ISSN={0003-9527},
     journal={Arch. Ration. Mech. Anal.},
      volume={243},
      number={3},
       pages={1151\ndash 1180},
      review={\MR{4381138}},
}

\bib{ELM23}{article}{
      author={{Elgindi}, Tarek~M.},
      author={{Liss}, Kyle},
      author={{Mattingly}, Jonathan~C.},
       title={{Optimal enhanced dissipation and mixing for a time-periodic,
  Lipschitz velocity field on $\mathbb{T}^2$}},
        date={2023-04},
     journal={arXiv e-prints},
      eprint={2304.05374},
}

\bib{ElgindiZlatos19}{article}{
      author={Elgindi, Tarek~M.},
      author={Zlato\v{s}, Andrej},
       title={Universal mixers in all dimensions},
        date={2019},
        ISSN={0001-8708},
     journal={Adv. Math.},
      volume={356},
       pages={106807, 33},
      review={\MR{4008523}},
}

\bib{FengIyer19}{article}{
      author={Feng, Yuanyuan},
      author={Iyer, Gautam},
       title={Dissipation enhancement by mixing},
        date={2019},
        ISSN={0951-7715},
     journal={Nonlinearity},
      volume={32},
      number={5},
       pages={1810\ndash 1851},
      review={\MR{3942601}},
}

\bib{IyerKiselevEA14}{article}{
      author={Iyer, Gautam},
      author={Kiselev, Alexander},
      author={Xu, Xiaoqian},
       title={Lower bounds on the mix norm of passive scalars advected by
  incompressible enstrophy-constrained flows},
        date={2014},
        ISSN={0951-7715},
     journal={Nonlinearity},
      volume={27},
      number={5},
       pages={973\ndash 985},
      review={\MR{3207161}},
}

\bib{IyerXuEA21}{article}{
      author={Iyer, Gautam},
      author={Xu, Xiaoqian},
      author={Zlato\v{s}, Andrej},
       title={Convection-induced singularity suppression in the
  {K}eller-{S}egel and other non-linear {PDE}s},
        date={2021},
        ISSN={0002-9947},
     journal={Trans. Amer. Math. Soc.},
      volume={374},
      number={9},
       pages={6039\ndash 6058},
      review={\MR{4302154}},
}

\bib{LunasinLinEA12}{article}{
      author={Lunasin, Evelyn},
      author={Lin, Zhi},
      author={Novikov, Alexei},
      author={Mazzucato, Anna},
      author={Doering, Charles~R.},
       title={Optimal mixing and optimal stirring for fixed energy, fixed
  power, or fixed palenstrophy flows},
        date={2012Nov.},
        ISSN={0022-2488},
     journal={J. Math. Phys.},
      volume={53},
      number={11},
}

\bib{LinThiffeaultEA11}{article}{
      author={Lin, Zhi},
      author={Thiffeault, Jean-Luc},
      author={Doering, Charles~R.},
       title={Optimal stirring strategies for passive scalar mixing},
        date={2011},
        ISSN={0022-1120},
     journal={J. Fluid Mech.},
      volume={675},
       pages={465\ndash 476},
      review={\MR{2801050 (2012h:76059)}},
}

\bib{HillSturmanEA22}{article}{
      author={Myers~Hill, J},
      author={Sturman, Rob},
      author={Wilson, Mark~CT},
       title={Exponential mixing by orthogonal non-monotonic shears},
        date={2022},
     journal={Physica D: Nonlinear Phenomena},
      volume={434},
       pages={133224},
}

\bib{pierrehumbert1994tracer}{article}{
      author={Pierrehumbert, RT},
       title={Tracer microstructure in the large-eddy dominated regime},
        date={1994},
     journal={Chaos, Solitons \& Fractals},
      volume={4},
      number={6},
       pages={1091\ndash 1110},
}

\bib{Thiffeault12}{article}{
      author={Thiffeault, Jean-Luc},
       title={Using multiscale norms to quantify mixing and transport},
        date={2012},
        ISSN={0951-7715},
     journal={Nonlinearity},
      volume={25},
      number={2},
       pages={R1\ndash R44},
      review={\MR{2876867}},
}

\bib{YaoZlatos17}{article}{
      author={Yao, Yao},
      author={Zlato\v{s}, Andrej},
       title={Mixing and un-mixing by incompressible flows},
        date={2017},
        ISSN={1435-9855},
     journal={J. Eur. Math. Soc. (JEMS)},
      volume={19},
      number={7},
       pages={1911\ndash 1948},
      review={\MR{3656475}},
}

\end{biblist}
\end{bibdiv}

\end{document}